\documentclass[reqno]{amsart}
\usepackage{fullpage}
\IfFileExists{stmaryrd.sty}{\usepackage{stmaryrd}}{}
\usepackage{graphicx}
\usepackage{amsmath,amssymb,amsfonts,mathrsfs,amscd}
\IfFileExists{newtxtext.sty}{%
  \IfFileExists{newtxmath.sty}{\usepackage{newtxtext,newtxmath}}{}%
}{}
\usepackage[all,cmtip]{xy}
\usepackage{hyperref}
\newtheorem{thm}{Theorem}[section]
\newtheorem{lem}[thm]{Lemma}

\theoremstyle{definition}

\newtheorem{exa}[thm]{Example}

\theoremstyle{remark}

\renewcommand{\baselinestretch}{1.3}\normalsize

\numberwithin{equation}{section}

\DeclareMathOperator{\Irr}{Irr}

\DeclareMathOperator{\Syl}{Syl}

\DeclareMathOperator{\Ker}{Ker}

\DeclareMathOperator{\Gal}{Gal}

\newtheorem*{questionD}{Question D}

\theoremstyle{plain}
\newtheorem{theoremA}{Theorem}

\newtheorem{theoremB}[theoremA]{Theorem}

\newtheorem{corollaryC}[theoremA]{Corollary}

\makeatletter
\renewcommand{\@setauthors}{%
  \begin{center}
    \normalfont\normalsize
    \def\and{,\space}%
    \authors
  \end{center}%
}
\makeatother

\makeatletter
\renewcommand{\section}{%
  \@startsection{section}{1}{\z@}%
    {.7\linespacing\@plus\linespacing}%
    {.5\linespacing}%
    {\normalfont\bfseries\raggedright
     \let\@secnumfont\bfseries}%
}
\makeatother

\begin{document}
\raggedbottom
\title{ON CHARACTER CONDUCTORS AND SYLOW
p-SUBGROUPS}
\author{Yang Zhou}
\author{Heguo Liu}

\maketitle

\noindent
\begin{minipage}{\textwidth}
\renewcommand{\baselinestretch}{1}
\normalfont\normalsize
\setlength{\parindent}{0pt}
\setlength{\parskip}{0pt}

\hrule height 0.4pt
\vspace{0.6em}

\textbf{Abstract}\par
\vspace{0.4em}

Let $P\in\Syl_p(G)$, where $G$ is a finite group and $p$ is a prime,
and let $\chi\in\Irr(G)$. In 2024, I. M. Isaacs and G. Navarro
asked whether the degree of $\mathbb{Q}_{p^f}$ over
$\mathbb{Q}_p(\chi_P)$ is always at most the $p$-part
$\chi(1)_p$ of $\chi(1)$, where $p^f$ is the $p$-part of the
conductor of $\chi$ and $f\geq1$. Here
$\mathbb{Q}_n=\mathbb{Q}(e^{2\pi i/n})$ is the $n$th cyclotomic
field, $\chi_P$ is the restriction of $\chi$ to $P$, and
$\mathbb{Q}_p(\chi_P)$ is the field generated over
$\mathbb{Q}_p$ by the values of $\chi_P$.
We prove that this inequality holds for $p$-solvable groups.
\par\medskip

\textit{Keywords:}\quad Character conductors, fields of values,
Sylow subgroups, $p$-solvable groups.
\par
\textit{MSC2020:}\quad 20C15, 20D20
\par\vspace{0.6em}

\hrule height 0.4pt

\hrule height 0.4pt
\end{minipage}
\par\medskip

\Large


\section{Introduction}

All groups considered in this paper are finite.
If $\chi$ is a character of a group $G$, we write
$\mathbb{Q}(\chi)$ to denote the field generated over $\mathbb{Q}$
by the values of $\chi$.
For a positive integer $n$, let
$\mathbb{Q}_n=\mathbb{Q}(e^{2\pi i/n})$ be the $n$th cyclotomic field.
The conductor $c(\chi)$ is the smallest positive integer $n$ such
that all values of $\chi$ lie in $\mathbb{Q}_n$.
For a prime $p$ and a positive integer $m$, we write $m_p$ to denote
the $p$-part of $m$. We say that $\chi$ has $p'$-degree if $p\nmid\chi(1)$.
If $H\leq G$, we write $\chi_H$ for the restriction of $\chi$ to $H$,
and $\mathbb{Q}_p(\chi_H)$ for the field generated over
$\mathbb{Q}_p$ by the values of $\chi_H$.

In 2019, I. M. Isaacs, M. W. Liebeck, G. Navarro, and P. H. Tiep
proved that if $\chi\in\Irr(G)$ has odd degree and
$\mathbb{Q}(\chi)=\mathbb{Q}(\sqrt{r})$, where $r\neq-1$ is
square-free, then $r\equiv1\pmod{4}$; see \cite{ILNT}.
This result was strengthened by G. Navarro and P. H. Tiep in 2021.
In particular, they proved that if $\chi$ has odd degree and
$c(\chi)_2=2^a$, then
$\mathbb{Q}_{2^a}\subseteq\mathbb{Q}(\chi)$; see \cite{NT21}.
These results show that there is a connection between the field of values of an irreducible character and the \(p\)-part of its degree.

Now fix a prime $p$ and let $P\in\Syl_p(G)$.
In 2024, Isaacs and Navarro proved that if $G$ is $p$-solvable and
$\chi\in\Irr(G)$ has $p'$-degree, and if $c(\chi)_p=p^f$,
where $f\geq1$, then
\[
  \mathbb{Q}_{p^f}=\mathbb{Q}_p(\chi_P);
\]
see \cite[Theorem~A]{IN}.
Since $\chi(1)_p=1$ in this situation, it is natural to ask the
following question, which appears as Question~D in \cite{IN}.

\begin{questionD}
Suppose that $G$ is a finite group, $p$ is a prime,
$P\in\Syl_p(G)$, and $\chi\in\Irr(G)$. Write
$c(\chi)_p=p^f$, where $f\geq1$. Is it true that
\[
\lvert\mathbb{Q}_{p^f}:\mathbb{Q}_p(\chi_P)\rvert
   \leq \chi(1)_p\,?
\]
\end{questionD}

We prove this inequality for $p$-solvable groups.

\begin{theoremA}\label{thm:A}
Let $\chi\in\Irr(G)$, where $p$ is a prime and $G$ is a
$p$-solvable group, and let $P\in\Syl_p(G)$. Write
$c(\chi)_p=p^f$, where $f\geq1$. Then
\[
 \left|\mathbb{Q}_{p^f}:\mathbb{Q}_p(\chi_P)\right|
 \leq \chi(1)_p.
\]
\end{theoremA}

This inequality has an application to the following
conjecture of N.~N.~Hung, G.~Navarro, and P.~H.~Tiep;
see \cite[Conjecture~B and Theorem~C]{HNT26}.
Let $F$ be an abelian extension of $\mathbb{Q}$ with conductor $c$,
and let $b$ be a nonnegative integer.
Then $F$ is the field of values of an irreducible character $\chi$
of some group with $\chi(1)_p=p^b$ if and only if
$|F\mathbb{Q}_{c_p}:F|_p\leq p^b$.
The necessity of this condition for $p$-solvable groups
was proved in Theorem~D(i) of \cite{HNT26}.
By the discussion following Question~5.8 in \cite{Hung},
Theorem~\ref{thm:A} gives another proof.

If $\chi(1)$ is a power of $p$, we obtain the following
stronger form of Theorem~\ref{thm:A}.

\begin{theoremB}\label{thm:B}
Let $\chi\in\Irr(G)$, where $p$ is a prime and $G$ is a finite
$p$-solvable group, and let $P\in\Syl_p(G)$. Suppose that
$\chi(1)$ is a power of $p$. Write $c(\chi)_p=p^f$, where
$f\geq1$. Then $\mathbb{Q}_p\subseteq\mathbb{Q}(\chi_P)$ and
\[
 \left|\mathbb{Q}_{p^f}:\mathbb{Q}(\chi_P)\right|
 [\chi_P,\chi_P]\leq\chi(1).
\]
\end{theoremB}

Theorem~\ref{thm:B} gives the following consequence for the case of
equality in Theorem~\ref{thm:A}.

\begin{corollaryC}\label{cor:C}
Let $\chi\in\Irr(G)$, where $G$ is a $p$-solvable group, and let
$P\in\Syl_p(G)$.
Suppose that $\chi(1)$ is a power of $p$, and write
$c(\chi)_p=p^f$, where $f\geq1$.
If
\[
  \left|\mathbb{Q}_{p^f}:\mathbb{Q}_p(\chi_P)\right|
  =\chi(1),
\]
then $\chi_P$ is irreducible.
\end{corollaryC}

For groups of odd order, there is a familiar condition on
fields of values that guarantees an irreducible restriction.
Suppose that $G$ has odd order and $\chi(1)$ is a power of $p$.
If $\mathbb{Q}(\chi)=\mathbb{Q}(\chi_P)$, then $\chi$ is
$p$-special, and hence $\chi_P$ is irreducible; see
\cite[Theorems~2.10 and~2.26]{I2018}. Corollary~\ref{cor:C} gives another
such condition: $\chi_P$ is irreducible when the index in
Theorem~\ref{thm:A} attains its upper bound.

\section{Selecting constituents}
\label{sec:selection}

We generalize Isaacs and Navarro's construction of the
$p$-group characters $\Delta_m(\Psi)$ in \cite{IN}.

Let $P$ be a $p$-group. For $\theta\in\Irr(P)$, let $m_\theta$
be the largest order of an element of $P/\Ker(\theta)$.
If $\Psi$ is a character of $P$ and $m$ is a power of $p$, define
\[
  \Delta_m(\Psi)=
  \sum_{\substack{\theta\in\Irr(P)\\m_\theta=m}}
  [\Psi,\theta]\theta.
\]
Thus $\Psi=\sum_m\Delta_m(\Psi)$, where $m$ runs over the distinct
values of $m_\theta$ for the irreducible constituents $\theta$ of $\Psi$.
If all constituents of $\Psi$ are linear, our construction agrees
with that of Isaacs and Navarro in \cite{IN}.

Now let $\chi\in\Irr(G)$, and let $P\in\Syl_p(G)$.
Write $\chi(1)_p=p^d$. We have
\[
  \chi(1)=\sum_m\Delta_m(\chi_P)(1).
\]
Since $p^{d+1}\nmid\chi(1)$, at least one term in this sum is
not divisible by $p^{d+1}$.
Let $n$ be the largest power of $p$ such that
\begin{equation}\label{eq:selection}
  p^{d+1}\nmid\Delta_n(\chi_P)(1).
\end{equation}

We first deduce Theorem~\ref{thm:A} from
Theorem~\ref{thm:selected-exponent}, whose proof is shown in
Section~\ref{sec:linear-rational}.

\begin{thm}\label{thm:selected-exponent}
Let $\chi\in\Irr(G)$, where $G$ is a $p$-solvable group,
and let $P\in\Syl_p(G)$.
Write $\chi(1)_p=p^d$ and $c=c(\chi)_p$, and let $n$ be the
largest power of $p$ such that
\[
  p^{d+1}\nmid\Delta_n(\chi_P)(1).
\]
If $c\nmid p^{d+1}$, then $c\mid n$.
\end{thm}

To prove Theorem~\ref{thm:A}, we need the following lemma.

\begin{lem}\label{lem:field-degree}
Let $\Psi$ be a character of a $p$-group $P$.
Suppose that $\Psi\neq0$ and that $m_\theta=m>1$ for every
irreducible constituent $\theta$ of $\Psi$.
If $\Psi$ is $F$-valued, where
$\mathbb{Q}_p\subseteq F\subseteq\mathbb{Q}_m$, then
$|\mathbb{Q}_m:F|\mid\Psi(1)$.
\end{lem}

\begin{proof}
Write $H=\Gal(\mathbb{Q}_m/F)$ and $e=|H|$.
Since $\mathbb{Q}_p\subseteq F$, the group $H$ is a $p$-group.
Also, $\Psi$ is $F$-valued, so $H$ permutes its irreducible
constituents, preserving their multiplicities.
It  suffices to show that $e$ divides the degree of the sum of the characters in each $H$-orbit.

Let $\mathcal{O}$ be such an orbit, choose
$\theta\in\mathcal{O}$, and write
\[
  \Gamma=\sum_{\eta\in\mathcal{O}}\eta.
\]
Then $\Gamma$ is $F$-valued, and
$\Gamma(1)=|\mathcal{O}|\theta(1)$ is a power of $p$.

Choose $x\in P$ such that $x\Ker(\theta)$ has order $m$.
Since $x^m\in\Ker(\theta)$, every linear constituent of
$\theta_{\langle x\rangle}$ has order dividing $m$.
Suppose that no constituent of $\theta_{\langle x\rangle}$
has order $m$. Since their orders divide $m$ and $m$ is
a power of $p$, they all divide $m/p$.
Thus $x^{m/p}$ lies in the kernel of every constituent,
and hence $x^{m/p}\in\Ker(\theta)$, contrary to the choice of $x$.

Choose a linear constituent $\lambda$ of
$\theta_{\langle x\rangle}$ of order $m$. The values of $\lambda$ generate $\mathbb{Q}_m$, so its
$H$-orbit has size $e$.
Since $\Gamma_{\langle x\rangle}$ is $F$-valued and contains
$\lambda$, all members of this orbit occur in
$\Gamma_{\langle x\rangle}$ with the same multiplicity.
These characters are linear, so $e\leq\Gamma(1)$.
Both integers are powers of $p$, and hence $e\mid\Gamma(1)$.

The constituents in $\mathcal{O}$ occur in $\Psi$ with
the same multiplicity, so their contribution to $\Psi(1)$
is divisible by $e$.
Since this holds for every $H$-orbit, we have $e\mid\Psi(1)$,
as required.
\end{proof}

If $\Psi$ is a character of $P$, then $\Delta_m(\Psi)$ is
$\mathbb{Q}(\Psi)$-valued.
To see this, suppose that $\sigma$ lies in the Galois group
$\Gal(\mathbb{Q}_{|P|}/\mathbb{Q}(\Psi))$.
Then $\sigma$ permutes the irreducible constituents of $\Psi$
and preserves their multiplicities and kernels.
Thus $\sigma$ fixes $\Delta_m(\Psi)$.

\begin{proof}[Proof of Theorem~\ref{thm:A},
  assuming Theorem~\ref{thm:selected-exponent}]
Write $\chi(1)_p=p^d$ and $c=c(\chi)_p$.
By \cite[Lemma~2.5]{IN}, we have
$\mathbb{Q}(\chi_P)\subseteq\mathbb{Q}_c$.
Since $p\mid c$, it follows that
$\mathbb{Q}_p\subseteq\mathbb{Q}_p(\chi_P)
\subseteq\mathbb{Q}_c$.

If $c\mid p^{d+1}$, then
\[
  |\mathbb{Q}_c:\mathbb{Q}_p(\chi_P)|
  \mid |\mathbb{Q}_{p^{d+1}}:\mathbb{Q}_p|
  =p^d,
\]
as required.
We can thus assume that $c\nmid p^{d+1}$.
Choose $n$ as in Theorem~\ref{thm:selected-exponent}, and write
$\Psi=\Delta_n(\chi_P)$.
Then $c\mid n$, and hence
$\mathbb{Q}_p(\chi_P)\subseteq\mathbb{Q}_n$.

By the choice of $n$, we have $\Psi\neq0$.
Also, $\Psi$ is $\mathbb{Q}_p(\chi_P)$-valued, and
$m_\theta=n$ for every irreducible constituent $\theta$ of $\Psi$.
It follows by Lemma~\ref{lem:field-degree} that
$|\mathbb{Q}_n:\mathbb{Q}_p(\chi_P)|\mid\Psi(1)$.
Since this index is a power of $p$ and
$p^{d+1}\nmid\Psi(1)$, we obtain
\[
  |\mathbb{Q}_c:\mathbb{Q}_p(\chi_P)|
  \mid |\mathbb{Q}_n:\mathbb{Q}_p(\chi_P)|
  \mid p^d.
\]
This proves the theorem.
\end{proof}

\section{Proof of Theorem~\ref{thm:selected-exponent}}
\label{sec:linear-rational}
In this section we prove Theorem~\ref{thm:selected-exponent}.
Before doing so, we make the following observation.

For a $p$-group $P$ and a power $r$ of $p$, write
$P^r=\langle x^r|x\in P\rangle$.
Let $\theta\in\Irr(P)$ and $N\lhd P$.
By Clifford's theorem, $[\theta_N,1_N]=\theta(1)$
if $N\subseteq\ker(\theta)$ and $[\theta_N,1_N]=0$ otherwise.
Since $m_\theta\mid r$ if and only if $P^r\subseteq\ker(\theta)$,
we have
\begin{equation}\label{eq:power-subgroup}
  [\Psi_{P^r},1_{P^r}]=\sum_{m\mid r}\Delta_m(\Psi)(1)
\end{equation}
for every character $\Psi$ of $P$.

Let $U\subseteq P$ and $L\lhd P$, and let $\eta$ be a
character of $U$.
By Mackey's formula, we have
$((1_L)^P)_U=|P:UL|(1_{U\cap L})^U$.
It follows by Frobenius reciprocity that
\begin{equation}\label{eq:mackey-inner-product}
\begin{aligned}
  \relax[(\eta^P)_L,1_L]
  &=[\eta^P,(1_L)^P]\\
  &=[\eta,((1_L)^P)_U]\\
  &=|P:UL|[\eta,(1_{U\cap L})^U]\\
  &=|P:UL|[\eta_{U\cap L},1_{U\cap L}].
\end{aligned}
\end{equation}

If $V\subseteq P$ and $h=|P:V|$, then $P^h\subseteq V$.
This follows by considering the action of $P$ on the cosets of $V$.
Now suppose that $s$ is a power of $p$ dividing both $r$
and $|P:U|$.
Choose $V$ with $U\subseteq V\subseteq P$ and $|P:V|=s$.
Then $P^r\subseteq P^s\subseteq V$, so $s\mid|P:UP^r|$.
By \eqref{eq:mackey-inner-product},
\begin{equation}\label{eq:index-divisibility}
  s\mid[(\eta^P)_{P^r},1_{P^r}].
\end{equation}

\begin{lem}\label{lem:linear-rational}
Let $P$ be a $p$-group, and let $U\subseteq P$.
Suppose that $\lambda$ is a linear character of $U$ and that
$\delta$ is a rational-valued character of $U$.
If $r$ is a $p$-power and $r<o(\lambda)$, then
$r\mid[((\lambda\delta)^P)_{P^r},1_{P^r}]$.
\end{lem}

\begin{proof}
Write $h=|P:U|$, $L=P^r$, and $\ell=o(\lambda)$.
If $r\leq h$, the result follows from \eqref{eq:index-divisibility}
with $s=r$.
We can thus assume that $h<r$.
Since $P^h\subseteq U$, we have
$U^r\subseteq L\subseteq U^{r/h}\subseteq U$.

Write $\tau=\sum_{\theta\in\Irr(U\mid\lambda_L)}
[\delta,\theta]\theta$.
Since $\lambda_L$ is $U$-invariant, it follows by Clifford's
theorem that $\theta_L=\theta(1)\lambda_L$ for each
constituent $\theta$ of $\tau$.
As $\delta$ is real valued, we have by
\eqref{eq:mackey-inner-product}
\[
  [((\lambda\delta)^P)_L,1_L]
  =h[\delta_L,\lambda_L]=h\tau(1).
\]
It therefore suffices to show that $r/h\mid\tau(1)$,
and we may assume that $\tau\neq0$.

Let $\theta$ be an irreducible constituent of $\tau$.
Since $r$ is a proper divisor of $\ell$ and both are powers of $p$,
we have $U^\ell\subseteq U^{\ell/p}\subseteq U^r\subseteq L$.
Since $\lambda$ is trivial on $U^\ell$ and
$\theta_L=\theta(1)\lambda_L$, we have
$U^\ell\subseteq\ker(\theta)$.
Thus every element of $U/\ker(\theta)$ has order dividing $\ell$.

Choose $u\in U$ such that $\lambda(u)$ has order $\ell$.
Since $u^{\ell/p}\in L$, we have
\[
  \theta(u^{\ell/p})
  =\theta(1)\lambda(u)^{\ell/p}
  \neq\theta(1).
\]
It follows that $u^{\ell/p}\notin\ker(\theta)$, and hence
$u\ker(\theta)$ has order $\ell$.
Thus $m_\theta=\ell$.

Since $L\subseteq U^{r/h}$, the character $\lambda_L$ is
$\mathbb{Q}_{\ell h/r}$-valued.
Every element of
$\Gal(\mathbb{Q}_{|U|}/\mathbb{Q}_{\ell h/r})$
fixes both $\delta$ and $\lambda_L$.
It permutes the irreducible constituents of $\delta$
lying over $\lambda_L$ and preserves their multiplicities.
Thus $\tau$ is $\mathbb{Q}_{\ell h/r}$-valued.
Since $p\mid \ell h/r$, it follows by
Lemma~\ref{lem:field-degree} that
\[
  \frac{r}{h}
  =|\mathbb{Q}_\ell:\mathbb{Q}_{\ell h/r}|
  \mid\tau(1).
\]
Hence $r\mid h\tau(1)$, as required.
\end{proof}

To apply Lemma~\ref{lem:linear-rational}, we use the following
decomposition of $\chi_P$.
Write $H^g=g^{-1}Hg$.

\begin{lem}\label{lem:sylow-decomposition}
Let $\chi\in\Irr(G)$, where $G$ is $p$-solvable, and let
$P$ be a Sylow $p$-subgroup of $G$.
Then there are a subgroup $T\subseteq P$, a linear character
$\lambda$ of $T$, and elements $g_1,\ldots,g_t\in G$ such that
\begin{equation}\label{eq:double-mackey}
  \chi_P=\sum_{i=1}^t(\lambda_i\delta_i)^P,
\end{equation}
where $U_i=T^{g_i}\cap P$,
$\lambda_i=(\lambda^{g_i})_{U_i}$, and $\delta_i$ is a
rational-valued character of $U_i$.
\end{lem}

\begin{proof}
Write $\chi=\gamma^G$, where $\gamma$ is a primitive character of a
subgroup $W$ of $G$. By \cite[Theorem~2.9]{I2018}, we can write
$\gamma=\alpha\beta$, where $\alpha$ is $p$-special and $\beta$ is
$p'$-special. Replacing $W$ and $\gamma$ by suitable conjugates,
we may assume that $Q=W\cap P$ is a Sylow $p$-subgroup of $W$.
Then $\alpha_Q$ is irreducible by \cite[Theorem~2.10]{I2018}.
Since irreducible characters of $p$-groups are monomial
(see \cite[Corollary~6.14]{I1976}), we can write
$\alpha_Q=\lambda^Q$, where $\lambda$ is a linear character of a
subgroup $T\subseteq Q$. Also, $\beta_Q$ is rational valued by
\cite[Corollary~2.14]{I2018}. Set $\eta=\lambda\beta_T$. Then
\[
  \gamma_Q=\alpha_Q\beta_Q=\lambda^Q\beta_Q
          =(\lambda\beta_T)^Q=\eta^Q.
\]

By Mackey's formula,
\[
  \chi_P=(\gamma^G)_P
        =\sum_x\bigl((\gamma^x)_{W^x\cap P}\bigr)^P,
\]
where $x$ runs over representatives for the double cosets $WxP$ in $G$.

We now consider one of these summands. Since $W\cap xPx^{-1}$
is a $p$-subgroup of $W$, we may choose $x$ in its double coset so that
$W\cap xPx^{-1}\subseteq Q$. Write $R=W^x\cap P$, so that
$R\subseteq Q^x$. Now $\gamma_Q=\eta^Q$, and hence $(\gamma^x)_R$
is the restriction of $(\eta^x)^{Q^x}$ to $R$. Applying Mackey's
formula again in $Q^x$, we deduce that this restriction is a sum of
characters
$\bigl((\eta^g)_U\bigr)^R$, where $g=xy$ for some $y\in Q^x$ and
$U=T^g\cap R$. Since $T^g\subseteq Q^x\subseteq W^x$, we have
$U=T^g\cap P$.

Inducing to $P$, we see that the summand
$\bigl((\gamma^x)_R\bigr)^P$ is a sum of characters of the form
\[
  \bigl((\eta^g)_U\bigr)^P
  =\bigl((\lambda^g)_U\,\delta\bigr)^P,
  \qquad \delta=(\beta^g)_U.
\]
Here $\delta$ is rational valued since it is a restriction of
$(\beta_T)^g$. Since this holds for every $x$, we have \eqref{eq:double-mackey} .
\end{proof}

Now we can prove Theorem~\ref{thm:selected-exponent}.

\begin{proof}[Proof of Theorem~\ref{thm:selected-exponent}]
Write $r=c/p$ and $L=P^r$. Since $c$ is a power of $p$ and
$c\nmid p^{d+1}$, we have $p^{d+2}\mid c$, and hence
$p^{d+1}\mid r$. Let $\Lambda=\sum_{m<c}\Delta_m(\chi_P)$.
Since $m$ and $c$ are powers of $p$, we have $m<c$ if and only if
$m\leq c/p=r$, or, equivalently, $m\mid r$.
By \eqref{eq:power-subgroup}, we have
\[
  \Lambda(1)=\sum_{m\mid r}\Delta_m(\chi_P)(1)=[\chi_L,1_L].
\]
We show that $p^{d+1}\mid\Lambda(1)$.

Choose $W,\alpha,\beta,Q,T$ and $\lambda$ as in the proof of
Lemma~\ref{lem:sylow-decomposition}, and write $\ell=o(\lambda)$.
Since $\chi=(\alpha\beta)^G$ and $\alpha_Q=\lambda^Q$, we have
\[
  \chi(1)=|G:W|\alpha(1)\beta(1)
  \qquad\text{and}\qquad
  \alpha(1)=|Q:T|.
\]
Now $\beta(1)$ is a $p'$-number and $|G:W|_p=|P:Q|$.
Taking $p$-parts, we obtain
\[
  p^d=|G:W|_p\,|Q:T|=|P:Q|\,|Q:T|=|P:T|.
\]

Next we show that $c\mid\ell$. If
$\sigma\in\Gal(\mathbb{Q}_{|W|}/\mathbb{Q}(\alpha_Q))$, then
$(\alpha^\sigma)_Q=(\alpha_Q)^\sigma=\alpha_Q$.
Since $\alpha^\sigma$ is also $p$-special, it follows by
Theorem~2.10 of \cite{I2018} that $\alpha^\sigma=\alpha$.
Thus $\mathbb{Q}(\alpha)\subseteq\mathbb{Q}(\alpha_Q)$, and the
reverse inclusion follows by restriction. The values of $\lambda$ lie in $\mathbb{Q}_\ell$.
Since $\alpha_Q=\lambda^Q$, we have
\[
  \mathbb{Q}(\alpha)=\mathbb{Q}(\alpha_Q)
  \subseteq\mathbb{Q}_\ell.
\]

Thus $c(\alpha)\mid\ell$ by Lemma~2.20 of \cite{I2018}.
Also, $c(\beta)$ is a $p'$-number by Corollary~2.13 of
\cite{I2018}. Since $\chi=(\alpha\beta)^G$, we have
\[
  c(\chi)\mid c(\alpha\beta)\mid c(\alpha)c(\beta),
\]
and we obtain $c\mid\ell$ by taking $p$-parts.

By Lemma~\ref{lem:sylow-decomposition}, $\chi_P$ is a sum of
characters $(\mu\delta)^P$, where $U=T^g\cap P$ for some
$g\in G$, $\mu=(\lambda^g)_U$, and $\delta$ is a rational-valued
character of $U$. It suffices to show that
\[
  p^{d+1}\mid[((\mu\delta)^P)_L,1_L]
\]
for each such summand.

Since $|T^g|=|T|$ and $|P:T|=p^d$, we have
$|P:U|=\frac{|P|}{|U|}
=\frac{|P|}{|T|}\frac{|T^g|}{|U|}
=|P:T|\,|T^g:U|=p^d\,|T^g:U|$.
If $|P:U|>p^d$, then $p^{d+1}\mid|P:U|$.
Since $p^{d+1}\mid r$, we can take $s=p^{d+1}$ in
\eqref{eq:index-divisibility} to obtain
$p^{d+1}\mid[((\mu\delta)^P)_L,1_L]$.

Suppose  that $|P:U|=p^d$. Then $|T^g:U|=1$, so
$U=T^g$ and $\mu=\lambda^g$ has order $\ell$.
Since $r=c/p$ and $c\mid\ell$, we see that $r$ is a proper divisor
of $\ell$. As $\delta$ is rational valued, we have
$r\mid[((\mu\delta)^P)_L,1_L]$ by
Lemma~\ref{lem:linear-rational}, and hence
$p^{d+1}\mid[((\mu\delta)^P)_L,1_L]$ in this case too.
Thus $p^{d+1}\mid\Lambda(1)$.

Finally,
\[
  \chi(1)=\Lambda(1)+\sum_{m\geq c}\Delta_m(\chi_P)(1).
\]
Since $\chi(1)_p=p^d$, it follows that
$p^{d+1}\nmid\Delta_m(\chi_P)(1)$ for some $m\geq c$.
By the choice of $n$, we have $n\geq m\geq c$, and thus
$c\mid n$, as required.
\end{proof}

\section{Characters of \texorpdfstring{$p$}{p}-power degree}
\label{sec:p-power-degree}

Let $\chi\in\Irr(G)$, where $G$ is any finite group, and let
$P\in\Syl_p(G)$. Write $c=c(\chi)_p$, and suppose that $c>1$.
Choose $x\in P$ such that $x\Ker(\chi)$ has the largest possible
order $m$. Then $c\mid m$, and $\chi_{\langle x\rangle}$
has a linear constituent of order $m$.
All Galois conjugates of this constituent over
$\mathbb{Q}_p(\chi_P)$ occur in $\chi_{\langle x\rangle}$.
Since these characters are linear, we obtain
\[
  |\mathbb{Q}_c:\mathbb{Q}_p(\chi_P)|
  \leq |\mathbb{Q}_m:\mathbb{Q}_p(\chi_P)|
  \leq\chi(1).
\]
Thus the inequality in Theorem~\ref{thm:A} holds for every
finite group when $\chi(1)$ is a power of $p$.

For $p$-solvable groups, we can say more about characters
of $p$-power degree. In this section we prove
Theorem~\ref{thm:B}, which strengthens this bound by including
the factor $[\chi_P,\chi_P]$.

\begin{proof}[Proof of Theorem~\ref{thm:B}]
We first show that $\mathbb{Q}_p\subseteq\mathbb{Q}(\chi_P)$.
Choose $W,\alpha,\beta,Q,T$ and $\lambda$ as in the proof of
Lemma~\ref{lem:sylow-decomposition}.
Since $\chi(1)=|G:W|\alpha(1)\beta(1)$ is a power of $p$,
we have $\beta(1)=1$ and $|G:W|=|P:Q|$.
Hence $G=WP$.
Also, $\beta_Q=1_Q$ since $\beta$ has $p'$-order.
By Mackey's formula, we have
\[
  \chi_P=(\alpha_Q)^P=(\lambda^Q)^P=\lambda^P.
\]
Write $c=c(\chi)_p$ and $\ell=o(\lambda)$.
As in the proof of Theorem~\ref{thm:selected-exponent},
we have $c\mid\ell$.
Since $c>1$, we have $\lambda\neq1_T$.
By Frobenius reciprocity,
\[
  [\chi_P,1_P]=[\lambda,1_T]=0.
\]

Put $F=\mathbb{Q}(\chi_P)$ and $e=|F\mathbb{Q}_p:F|$.
Let $\theta$ be an irreducible constituent of $\chi_P$.
Since $\theta\neq1_P$, we can choose $u\in P$ such that
$u\ker(\theta)$ is a central element of order $p$ in
$P/\ker(\theta)$.
By Schur's lemma, this element acts as a scalar of order $p$
in a faithful representation of $P/\ker(\theta)$
affording $\theta$.
Thus $\theta(u)/\theta(1)$ is a primitive $p$-th root of unity,
and hence $\mathbb{Q}_p\subseteq\mathbb{Q}(\theta)$.

Write $E=\mathbb{Q}_{|P|}$ and $H=\Gal(E/F)$.
An element of $H$ fixes $\theta$ if and only if it fixes
all values of $\theta$.
Thus the stabilizer of $\theta$ in $H$ is $\Gal(E/F(\theta))$.
Here $F(\theta)=F\mathbb{Q}(\theta)$, so
$F\mathbb{Q}_p\subseteq F(\theta)$.
The $H$-orbit of $\theta$ has size
\[
  |H:\Gal(E/F(\theta))|
  =|F(\theta):F|
  =e\,|F(\theta):F\mathbb{Q}_p|.
\]
Since $\chi_P$ is $F$-valued, $H$ permutes its irreducible
constituents and preserves their degrees and multiplicities.
Every such orbit has size divisible by $e$, so $e\mid\chi(1)$.
But $e\mid p-1$ and $\chi(1)$ is a power of $p$,
so $e=1$.
Hence $\mathbb{Q}_p\subseteq F$.

We now prove the inequality.
By Lemma~2.5 of \cite{IN}, we have $F\subseteq\mathbb{Q}_c$.
Since $c\mid\ell$, we also have
$\mathbb{Q}_c\subseteq\mathbb{Q}_\ell$.
Put $b=[\chi_P,\chi_P]$.
Since $\chi_P=\lambda^P$, we have by Frobenius reciprocity
\[
  b=[\chi_P,\chi_P]=[\chi_P,\lambda^P]=[\chi_T,\lambda].
\]
Thus $\lambda$ occurs in $\chi_T$ with multiplicity $b$.
The values of $\lambda$ generate $\mathbb{Q}_\ell$, so its
stabilizer in $H$ is $\Gal(E/\mathbb{Q}_\ell)$.
Its $H$-orbit has size $|\mathbb{Q}_\ell:F|$.
Since $\chi_T$ is $F$-valued, every member of this orbit
occurs in $\chi_T$ with multiplicity $b$.
All these characters are linear, so
\[
  |\mathbb{Q}_c:F|[\chi_P,\chi_P]
  \leq |\mathbb{Q}_\ell:F|\,b
  \leq\chi(1),
\]
as required.
\end{proof}

\begin{proof}[Proof of Corollary~\ref{cor:C}]
By Theorem~\ref{thm:B},
$\mathbb{Q}_p(\chi_P)=\mathbb{Q}(\chi_P)$, and hence
\[
  \chi(1)[\chi_P,\chi_P]
  =|\mathbb{Q}_{p^f}:\mathbb{Q}(\chi_P)|[\chi_P,\chi_P]
  \leq\chi(1).
\]
Thus $[\chi_P,\chi_P]=1$, and $\chi_P$ is irreducible.
\end{proof}
We give an example showing that the equality assumed in
Corollary~\ref{cor:C} can occur.

\begin{exa}\label{ex:corollary-C}
Let $G$ be the dihedral group of order $24$, let
$P\in\Syl_2(G)$, and let $\chi\in\Irr(G)$ be faithful
of degree $2$. Then
\[
  \mathbb{Q}(\chi)=\mathbb{Q}(\sqrt{3}),\qquad
  c(\chi)=12,\qquad
  \mathbb{Q}(\chi_P)=\mathbb{Q}.
\]
Thus
\[
  |\mathbb{Q}_{c(\chi)_2}:\mathbb{Q}_2(\chi_P)|
  =|\mathbb{Q}_4:\mathbb{Q}|
  =2=\chi(1).
\]
\end{exa}

\section*{Acknowledgements}

The authors thank Junwei Zhang and Zhenzhao Huang for
helpful suggestions.

At an early stage of this work, the authors used ChatGPT
(OpenAI) in an attempt to prove Theorem~\ref{thm:A}.
It found a proof only in the case where $P$ is abelian.
This argument suggested to the authors how to generalize
the construction of I.~M. Isaacs and G. Navarro in \cite{IN}.
ChatGPT was also used to find Example~\ref{ex:corollary-C}
and to improve the wording of some proofs.
The authors developed the proof for the general
$p$-solvable case and the remaining mathematical arguments
in this paper. They take full responsibility for the
content of this paper.

The authors acknowledge support from the National Natural Science Foundation of China (12571017).



School of Mathematics and Statistics, Hainan University, Haikou 570228, China\par
Email address: zhouyang2026@hainanu.edu.cn (Yang Zhou)\par

\medskip

School of Mathematics and Statistics, Hainan University, Haikou 570228, China\par
Email address: ghliu@hainanu.edu.cn (Heguo Liu)\par

\end{document}